\newif\ifpdf\ifx\pdfoutput\undefined\pdffalse%
                \documentclass{article}\usepackage{graphicx}%
        \else\pdfoutput=1\pdftrue%
                \documentclass[pdftex]{article}%
                \usepackage[pdftex]{hyperref}%
                \usepackage[pdftex]{graphicx}%
        \fi

\usepackage{color}
\usepackage{amsmath}
\usepackage{amsfonts}
\usepackage{amssymb}
\usepackage[normalem]{ulem}
\usepackage[latin1]{inputenc}

\newtheorem{theorem}{Theorem}

\newtheorem{corollary}[theorem]{Corollary}

\newenvironment{proof}[1][Proof]{\textbf{#1.} }%
               {\ \hfill\rule{0.5em}{0.5em}\medskip}

\newcommand{\defby}{\stackrel{{\scriptscriptstyle\text{def}}}{=}}
\newcommand{\PP}{\mathbb{P}}
\newcommand{\EE}{\mathbb{E}}
\newcommand{\ind}{\mathbf{1}}
\newcommand{\Psiper}{\Psi_{\textrm{per}}}
\newcommand{\tildeDelta}{\tilde\Delta}

\begin{document}
\title{Pinning by a sparse potential}
\author{
É. Janvresse, T. de la Rue, Y. Velenik\\
Laboratoire de Mathématiques Raphaël Salem\\
UMR-CNRS 6085, Université de Rouen\\
\small\texttt{Elise.Janvresse@univ-rouen.fr,} \\
\small\texttt{Thierry.de-la-Rue@univ-rouen.fr,} \\
\small\texttt{Yvan.Velenik@univ-rouen.fr}
}
\maketitle

\begin{abstract}
We consider a directed polymer interacting with a diluted pinning potential restricted to a line. We characterize explicitely the set of disorder configurations that give rise to localization of the polymer. We study both relevant cases of dimension $1+1$ and $1+2$. We also discuss the case of massless effective interface models in dimension $2+1$.
\end{abstract}




\bigskip

It is customary, when modelling a disordered physical system, to assume that the disorder is sampled from some suitable random distribution. Of course there is a high degree of arbitrariness in the choice of this distribution, and one hopes that only qualitative features are relevant. Then, in the best possible cases, one can prove results that hold for almost every disorder configuration. However, there are several drawbacks with such an approach : First, it would be desirable to avoid these additional assumptions on the distribution of the disorder, and second, even with an almost sure result, we are left clueless about the validity of the desired property when an explicit disorder configuration is given. Therefore, it would be very valuable if one could instead characterize the set of realizations of the environment for which a specific property holds, or at least give some sufficient conditions. This is a much more ambitious program, and it is probably doomed to fail in general. In this paper, we give a simple example of a problem where such an approach can actually be pursued.

An important physical problem, which has received much attention recently from the mathematical community, is that of a directed polymer in a random environment, see, e.g., \cite{CoShYo04} and references therein. The latter is modelled by an exponential perturbation of the path measure of a $d$-dimensional random walk (or Brownian motion), depending on the realization of a random environment. The perturbation is such that visits of the random walk to regions where the environment takes positive values are rewarded, while visits to regions of opposite sign are penalized. Among the many questions of interest, there is the problem of studying the superdiffusive behaviour of the path in dimension $1+1$\footnote{We use the terminology ``dimension $d_1+d_2$" when considering a $d_1$-dimensional (directed) object in a $(d_1+d_2)$-dimensional space.}. Heuristically, one expects that there will be an ``optimal tube'' in the environment, inside which the landscape looks particularly good from the random walk point of view, and along which the random walk will localize. It is natural to split this problem into two : 1) Establish the existence of such an ``optimal tube'', whatever that really means, and 2) Prove that given such an ``optimal tube'', there is pinning of the polymer along the tube. Once we accept that such a splitting is natural, one can start to build simpler models for both points separately, in order to gain a better understanding of these issues. For the second part, a natural simplification is as follows : Consider a path $Y : \mathbb{N}\to\mathbb{Z}^d$, and perturb the path measure of a random walk $X$ by rewarding each intersection between the paths $X$ and $Y$. The question is then to understand under which conditions on the path $Y$ there will be pinning of the polymer $X$, i.e, there will be a positive density of such intersections. This might then shed some light on the properties of this ``optimal tube'' one should look for when analyzing the random environment. The only result we are aware of in this direction is due to Ioffe and Louidor~\cite{IoLo04}, who consider the situation where $Y$ is itself a random walk (its increments having possibly a different law from those of $X$), and proved that pinning occurs in dimension $1+1$ for almost all realizations of $Y$. This is however ``only'' an almost sure result, and as such it does not tell what is the set of measure $0$ (w.r.t. the law of the random walk $Y$) of paths which do not lead to pinning, which is most unfortunate since the ``optimal tube'' is expected to behave quite differently from a random walk trajectory. It would thus be very interesting to get an explicit characterization of the paths $Y$ for which pinning occurs, or at least sufficient conditions.

We hope to come back to this issue in the future. In the present work, we analyze a much simpler situation. Namely, here the pinning potential is restricted to a single line, and the disorder comes from the fact that this potential is diluted. More precisely, let $\omega\in\{0,1\}^\mathbb{N}$, $\eta>0$, and let $\PP_0$ denote the law of an aperiodic, symmetric random-walk on $\mathbb{Z}^d$ starting from $0$, with increments of finite variance. Our main interest is in the following perturbation of $\PP_0$,
\begin{equation}
\PP_{N,\eta}^\omega (X) \defby (Z_{N,\eta}^\omega)^{-1} \;\exp\left(\eta\sum_{i\in\Lambda_N} \ind_{(X_i=0)} \omega_i \right) \, \PP_0 (X), 
\end{equation}
where $\Lambda_N\defby \{1,\ldots, N\}$ and $Z_{N,\eta}^\omega$ is the partition function used to normalize $\PP_{N,\eta}^\omega$ to a probability measure. This 
measure models the interaction of a directed polymer in $1+d_2$ dimensions, interacting with an attractive diluted potential restricted to the line $x_2=\dots=x_{d_2+1} = 0$. The central question is under which conditions does such a potential localize the polymer, i.e., when is it true that
$$
\liminf_{N\to\infty} \frac1N \EE_{N,\eta}^\omega\left( \sum_{i=1}^N \ind_{X_i = 0} \right) > 0 \text{ ?}
$$
When this happens, we say that there is pinning of the polymer by the potential.
It has been known for a long-time that in the special case $\omega\equiv 1$, the polymer is pinned for any value of $\eta>0$ when $d_2=1$ or $d_2=2$, but is not pinned for small enough values of $\eta$ in higher dimensions (see \textit{e.g.}~\cite{Bu81} for a special case, and \cite{BoVe01} for a more general treatment). The general case of a discrete-time Markov chain interacting with a (possibly
random) potential restricted to a line was recently investigated by Alexander
and Sidoravicius \cite{AlSi05}. In particular, they compared the effect of an i.i.d. random
potential with the constant potential given by its average. One of their main
results is that pinning of the polymer is strictly enhanced by the presence of
such disorder. The situation we consider is a simple particular case of their setting, but our result is stronger since
we work with a fixed (arbitrary) environment.
 
Of course, diluting the potential only makes it less likely for the polymer to be pinned, so there is still delocalization at small values of $\eta$ in dimensions $3$ and higher, for arbitrary $\omega$. In this work, we therefore restrict our attention to dimensions $1$ and $2$. Rather remarkably, in this case it is possible to obtain a very simple characterization of the set of environments for which pinning occurs, see Theorem~\ref{thm_Main} below and its corollary. Before stating the result, we also introduce another case where the same question can be investigated: $2+1$-dimensional massless effective interface models. In this case, let $\Lambda\Subset\mathbb{Z}^2$, let $V:\mathbb{R}\to\mathbb{R}$ such that $0 < c_- \leq V'' \leq c_+ < \infty$, and let $\eta\geq 0$. We are interested in the measure (on $\mathbb{R}^{\Lambda}$) defined by\footnote{In this paper, we only consider the so-called $\delta$-pinning. However, given the very rough nature of the bounds we are after, there is no difficulty in treating more general potentials.}
\begin{multline}
\label{EffInt}
\PP_{\Lambda,\eta}^\omega (dX) \defby (Z_{\Lambda,\eta}^\omega)^{-1} \;\exp \Bigl( -\tfrac12\sum_{\substack{i \sim j\\ i,j\in\Lambda}} V(X_i-X_j) + \sum_{\substack{i \sim j\\ i\in\Lambda,j\not\in\Lambda}} V(X_i) \Bigr) \\
\times \prod_{i\in\Lambda} \left( dX_i + \eta \omega_i \delta_0(dX_i) \right) ,
\end{multline}
where $\delta_0$ is the Dirac mass at $0$ and $\omega\in\{0,1\}^{\mathbb{N}^2}$. In the special case $\Lambda = \Lambda_N \defby \{1,\ldots,N\}^2$, we simply write $\PP_{N,\eta}^\omega$ and $Z_{N,\eta}^\omega$. This models a two-dimensional interface in a three-dimensional medium, interacting with an attractive (diluted) potential located in the plane $x_3=0$. The basic question is the same as above: Determine under which conditions the interface is localized by the potential. The case $\omega\equiv1$ has been studied in details recently, see~\cite{BoVe01} and reference therein. It turns out that in this case too, an explicit description of the set of disorder configurations leading to pinning can be obtained.

\begin{figure}
	\centering
	\includegraphics[width=12cm,height=2cm,bb=0 0 613 209]{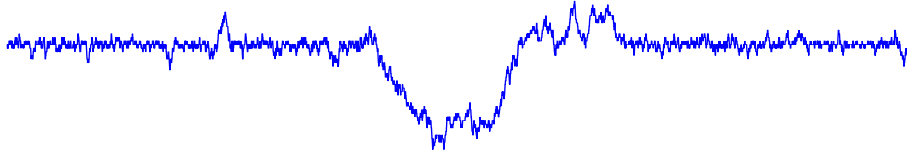}
	\caption{A simulation of the 1+1-dimensional process. The environment has density $.8$ on the first and last third of the interval, and $0$ inbetween.}
	\label{fig:path}
\end{figure}

\medskip
The following theorem and its corollary are valid for all three cases of dimensions $1+1$, $1+2$ and $2+1$.
\begin{theorem}
\label{thm_Main}
Let $0<\delta<1$ and $\eta>0$. For all $N>N_0(\delta,\eta)$, and for all $\omega$ such that
\begin{equation}
\label{density}
\sum_{i\in\Lambda_N}\omega_i > \delta |\Lambda_N|,
\end{equation}
we have for some $C=C(\delta,\eta)>0$
\begin{equation}
\EE_{\eta,N}^\omega\left[\sum_{i\in\Lambda_N} \ind_{(X_i=0)}\,\omega_i\right] > C |\Lambda_N|.
\end{equation}
\end{theorem}
\begin{corollary}
\label{cor_Main}
For any $\eta>0$,
\begin{equation}
\label{eq_densityOfVists}
\liminf_{N\to\infty} |\Lambda_N|^{-1} \EE_{\eta,N}^\omega\left[\sum_{i\in\Lambda_N}^N\ind_{(X_i=0)}\,\omega_i\right] > 0 ,
\end{equation}
if and only if
\begin{equation}
\liminf_{N\to\infty} |\Lambda_N|^{-1} \sum_{i\in\Lambda_N} \omega_i > 0.
\end{equation}
\end{corollary}
Note that this is a sensible definition of pinning, since if pinning does not hold in this sense, then there exists a sequence of increasing boxes, such that along this sequence the density of pinned sites goes to zero. Of course, along other sequences there can be a strictly positive density of pinned sites (examples are easily constructed).

\medskip
\begin{proof}[Proof of Theorem~\ref{thm_Main}]
\textsc{Step 1.} It is enough to prove that $Z_{N,\eta}^\omega/Z_{N,0}^\omega>D^{|\Lambda_N|}$, for some $D=D(\delta,\eta)>1$.
Indeed,
\begin{equation}
\label{eq_identity}
\log \frac{Z_{N,\eta}^\omega}{Z_{N,0}^\omega} = \int_0^\eta \EE_{\tilde\eta,N}^\omega\left[\sum_{i\in\Lambda_N} \ind_{(X_i=0)}\,\omega_i\right]\,d\tilde\eta,
\end{equation}
and, since the expectation is increasing in $\tilde\eta$, we obtain
$$ \EE_{\eta,N}^\omega\left[\sum_{i\in\Lambda_N} \ind_{(X_i=0)}\,\omega_i\right] > \dfrac{|\Lambda_N|\log D}{\eta}. $$
Of course, in the cases of dimensions $1+1$ and $1+2$, $Z_{N,0}^\omega=1$.

\smallskip
\textsc{Step 2.} We first treat the cases of polymers, that is dimensions $1+d_2$, $d_2=1$ or $2$. Let 
$$ \Omega_N\defby\{i\in\Lambda_N,\ \omega_i=1\}.$$
Writing 
$$\exp\left(\eta\sum_{i\in\Lambda_N} \ind_{(X_i=0)}\, \omega_i\right)=\prod_{i\in\Lambda_N} \left((e^\eta-1)\ind_{(X_i=0)}\ind_{i\in\Omega_N}+1\right) $$
and expanding the product, we get
\begin{align}
\nonumber
Z_{N,\eta}^\omega & = \EE_0\left[\exp\left(\eta\sum_{i\in\Lambda_N} \ind_{(X_i=0)}\, \omega_i\right)\right]\\
& = \sum_{A\subset\Omega_N}(e^\eta-1)^{|A|}\,\PP_0(X_i\equiv 0\text{ on }A).
\label{decomposition}
\end{align}

\smallskip
\textsc{Step 2.1.}
It is convenient to number the sites of $\Omega_N$, $\Omega_N=\{t_1<\cdots<t_{|\Omega_N|}\}$. Restricting the sum to the subsets $A$ of fixed cardinality $r$ (to be chosen later), and using the Markov property
for the random walk, we obtain the following lower bound for the partition function (where $t_0=\ell_0\defby 0$). 
\begin{equation}
\label{lower-bound}
Z_{N,\eta}^\omega \ge (e^\eta-1)^r\sum_{0<l_1<\cdots<l_r\le|\Omega_N|}\,\prod_{i=1}^r\PP_0\left(X_{t_{\ell_i}-t_{\ell_{i-1}}}=0\right),
\end{equation}
which yields, by the local CLT, for some $c>0$,
\begin{equation}
\label{lower-bound2}
Z_{N,\eta}^\omega \ge (e^\eta-1)^r\sum_{0<l_1<\cdots<l_r\le|\Omega_N|}\,\prod_{i=1}^r\dfrac{c}{(t_{\ell_i}-t_{\ell_{i-1}})^{d_2/2}}.
\end{equation}

\smallskip
\textsc{Step 2.2.} We begin by considering the simpler case of dimension $1+1$. Observe that, by Jensen inequality, we have 
\begin{equation}
\label{jensen}
\prod_{i=1}^r\dfrac{1}{\sqrt{t_{\ell_i}-t_{\ell_{i-1}}}} = \exp\left(-\dfrac{r}{2}\sum_{i=1}^r\dfrac{1}{r}\log(t_{\ell_i}-t_{\ell_{i-1}})\right)\ge \exp\left(-\dfrac{r}{2}\log\dfrac{N}{r}\right).
\end{equation}
Setting $r=|\Omega_N|/K$ for some integer $K$ to be chosen later, and observing that the number of terms in the RHS of \eqref{lower-bound2}
is at least $K^r$, we obtain
$$   Z_{N,\eta}^\omega \ge \left(Kc(e^\eta-1)\sqrt{\dfrac{r}{N}}\right)^r. $$
Using $|\Omega_N|\ge \delta N$, and choosing $K$ large enough, the conclusion follows. 

\smallskip
\textsc{Step 2.3.} We now turn to the more delicate case of dimension $1+2$. Although the above argument involving Jensen inequality is too rough to conclude now, it suggests that the worst possible environment $\omega$ with a fixed density $\delta$
occurs when $t_i-t_{i-1}\equiv \delta^{-1}$. 

We introduce $\Delta_i \defby t_i-t_{i-1}$, and
\begin{align*}
\Psi(\Delta_1,\ldots,\Delta_{|\Omega_N|}) &\defby
\sum_{0<\ell_1<\cdots<\ell_r\le|\Omega_N|} \,\prod_{i=1}^r\dfrac{1}{(t_{\ell_i}-t_{\ell_{i-1}})}\\
& =\sum_{0<\ell_1<\cdots<\ell_r\le|\Omega_N|} \frac1{(\Delta_1+\cdots+\Delta_{\ell_1})\cdots (\Delta_{\ell_{r-1}+1}+\cdots+\Delta_{\ell_r})} .
\end{align*}
Instead of working directly with the function $\Psi$, it is convenient to consider a periodized version defined by
(see Figure~\ref{figure})
\begin{multline*}
\Psiper(\tildeDelta_1,\Delta_2,\ldots,\Delta_{|\Omega_N|}) \defby\\ \sum_{0<\ell_1<\cdots<\ell_r\le|\Omega_N|} \frac1{(\Delta_{\ell_r+1}+\cdots+\Delta_{|\Omega_N|}+\tildeDelta_1+\Delta_2+\cdots+\Delta_{\ell_1})}\times\\
\prod_{i=2}^r \frac1{(\Delta_{\ell_{i-1}+1}+\cdots+\Delta_{\ell_i})} ,
\end{multline*}
where $\tildeDelta_1 \defby N+1 - \sum_{i=2}^{|\Omega_N|}\Delta_i$.
\begin{figure}
\begin{center}\input{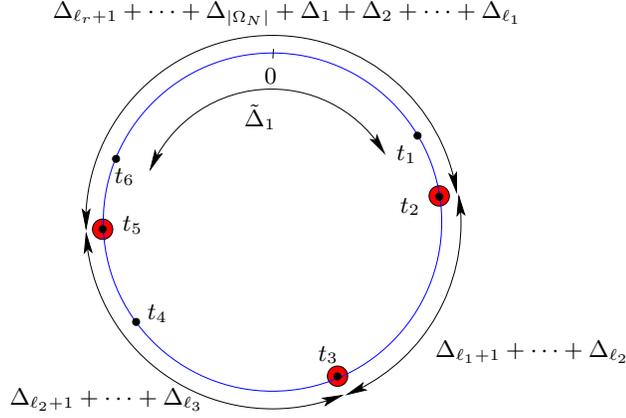}
\end{center}
\caption{The construction of $\Psiper$: Here $|\Omega_N|=6$, $r=3$, $\ell_1=2$, $\ell_2=3$, $\ell_3=\ell_r=5$.}
\label{figure}
\end{figure} 
We are going to determine the (unique) minimum of the function $\Psiper$, seen as a function on $\mathbb{R}_+^{|\Omega_N|}$, restricted to the manifold $\tildeDelta_1+\sum_{i=2}^{|\Omega_N|} \Delta_i = N+1$.

Notice first that $\Psiper$ is a convex function. Indeed, the function $(x_1,\ldots,x_k) \mapsto (x_1\cdots x_k)^{-1}$ is convex on $\mathbb{R}_+^k$, the composition of a convex function with an affine function is convex, and the sum of convex functions is also convex, as well as its restriction to an affine subspace.

We claim that the point $\Delta_i \equiv \Delta \defby (N+1)/|\Omega_N|$ is a (strict) local minimum of $\Psiper$, and therefore its unique minimum. To prove this, it is enough, by symmetry, to show that
\begin{equation}
\label{minloc}
\Psiper(\Delta,\Delta,\Delta,\ldots,\Delta) - \Psiper(\Delta+h,\Delta-h,\Delta,\ldots,\Delta) \leq 0 ,
\end{equation}
for all $h$ small enough. Indeed, each allowed configuration $(\tildeDelta_1,\Delta_2,\ldots,\Delta_{|\Omega_{N}|})$
can be written as $(\Delta+h_1,\Delta+h_2-h_1,\ldots,\Delta+h_{|\Omega_{N}|-1}-h_{|\Omega_{N}|-2},\Delta-h_{|\Omega_{N}|-1})$, and equation~\eqref{minloc} together with invariance of $\Psiper$ under cyclic permutation of the variables, ensure that all partial derivatives with respect to the $h_i$'s are nonnegative.

Observing that $\Psiper$ is also invariant under the transformation $$(\tildeDelta_1,\Delta_2,\ldots,\Delta_{|\Omega_N|}) \mapsto (\Delta_{|\Omega_N|},\ldots,\Delta_2,\tildeDelta_1), $$
we have that
$$
\Psiper(\Delta+h,\Delta-h,\Delta,\ldots,\Delta) = \Psiper(\Delta-h,\Delta+h,\Delta,\ldots,\Delta) .
$$
Therefore, the claim follows by convexity.

We need to compare $\Psi$ and $\Psiper$. Noticing that $\Delta_{\ell_r+1}+\cdots+\Delta_{|\Omega_N|}+\tildeDelta_1 \geq \Delta_1$, we immediately get that
$$
\Psi(\Delta_1,\ldots,\Delta_{|\Omega_N|}) \geq \Psiper(\tildeDelta_1,\Delta_2,\ldots,\Delta_{|\Omega_N|}) \geq \Psiper(\Delta,\ldots,\Delta) .
$$
Therefore $Z_{N,\eta}^\omega \geq \left( c(e^\eta -1) \right)^r \, \Psiper(\Delta,\ldots,\Delta)$. It only remains to find a bound on $\Psiper(\Delta,\ldots,\Delta)$. Let $K>0$; this number will be chosen later. We have that, for all $N$ large enough,
\begin{align*}
\Psiper(\Delta,\ldots,\Delta)
&\geq \frac1N \sum_{\substack{0<\ell_1<\ell_2<\ldots<\ell_r\\|\ell_i-\ell_{i-1}| \leq K}}\, \prod_{i=2}^r \frac1{(\ell_i-\ell_{i-1})\Delta}\\
&= \frac1{\Delta^{r-1}N} \sum_{\substack{0<\ell_1<\ell_2<\ldots<\ell_{r-1}\\|\ell_i-\ell_{i-1}| \leq K}}\, \prod_{i=2}^{r-1} \frac1{\ell_i-\ell_{i-1}} \, \sum_{k=1}^K \frac1k\\
&\geq \frac1{\Delta^{r-1}N} \log K \, \sum_{\substack{0<\ell_1<\ell_2<\ldots<\ell_{r-1}\\|\ell_i-\ell_{i-1}| \leq K}}\, \prod_{i=2}^{r-1} \frac1{\ell_i-\ell_{i-1}}\\
&\geq \frac{(\log K)^{r-1}}{\Delta^{r-1} N} .
\end{align*}
Choosing now $K$ large enough (which is possible as soon as $N>N_0(\delta,\eta)$), we conclude the proof of the theorem.

\smallskip

\begin{figure}
\label{figure2}
\begin{center}\input{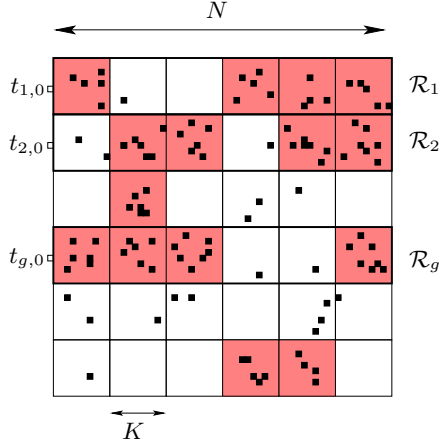}
\end{center}
\caption{The partition of $\Lambda_N$: Good cells are shaded, good rows are labelled $\mathcal{R}_1,\ldots,\mathcal{R}_g$. 
We also indicate the extra boundary points $t_{1,0},\ldots,t_{g,0}$.}
\end{figure} 

\textsc{Step 3.} We finally consider the case of dimension $2+1$. Let $\Omega_N \defby\{i\in\Lambda_N, \omega_i=1\}$. Expanding the product in~\eqref{EffInt}, we obtain a representation similar to~\eqref{decomposition},
\begin{equation}
\label{expansion}
\frac{Z_{N,\eta}^\omega}{Z_{N,0}^\omega} = \sum_{A\subset\Omega_N} \eta^{|A|}\, \frac {Z_{\Lambda_N \setminus A,0}^\omega}{Z_{N,0}^\omega} .
\end{equation}
Let us partition $\Lambda_N$ into cells of sidelength $K$ (to be chosen later). We suppose, to ease notations, that this partitionning can be done exactly; the general case is treated in a straightforward way.

Let $0<\rho<\delta/(2-\delta)$, and let us say that a cell is \textit{good} if it contains at least $\rho K^2$ sites of $\Omega_N$. Clearly, there is at least a density $\rho/(1+\rho)$ of good cells, since otherwise
$$
\sum_{i\in\Lambda_N} \omega_i \leq \left( 1-\frac{\rho}{1+\rho} \right) \frac{|\Lambda_N|}{K^2}\, \rho K^2 + \frac{\rho}{1+\rho}\frac{|\Lambda_N|}{K^2}\, K^2 = \frac{2\rho}{1+\rho} |\Lambda_N| <\delta |\Lambda_N| .
$$
Now, let us say that a row of cells is \textit{good} if the number of good cells in this row is at least $\zeta N/K$, where $\zeta$ is some small enough constant. A similar computation shows that for the class of environments we consider, there must be at least a fraction $\zeta/(1+\zeta)$ of good rows.

Returning to~\eqref{expansion}, we see that we must find a reasonable lower bound on the ratio of partition functions in the RHS, for a large enough class of sets $A$. Let us denote by $\mathbf{A}$ the class of sets $A$ containing exactly one site in each good cell located in a good row. The good rows can be numbered $\mathcal{R}_1,\ldots,\mathcal{R}_g$, with $g\geq \zeta/(1+\zeta) N/K$. $A_k$, the set of sites of $A\in\mathbf{A}$ belonging to the row $\mathcal{R}_k$, can then also be ordered according to their first coordinate, $A_k = \{t_{k,1},\ldots,t_{k,n_k}\}$, where $n_k\geq \zeta N/K$. For each $k$, let also $t_{k,0}$ be a site of $\mathbb{Z}^2\setminus\Lambda_N$ neighbour of the leftmost cell of $\mathcal{R}_k$. We need the following result from~\cite{DeVe00} : For any $B\Subset\mathbb{Z}^2$ and $t\in B$,
$$
\frac {Z_{B\setminus \{t\},0}^\omega}{Z_{B,0}^\omega} \geq \frac{c}{\sqrt{\log(1+d(t,B^c))}}.
$$
From this, we obtain, setting $A_{k,i} = A \setminus \{ t_{\ell,j} : \ell<k, \text{ or } \ell=k \text{ and } j\leq i \}$, that
\begin{align*}
\frac {Z_{\Lambda_N \setminus A,0}^\omega}{Z_{N,0}^\omega}
&= \prod_{k=1}^g \prod_{i=1}^{n_k} \frac {Z_{\Lambda_N \setminus A_{k,i-1},0}^\omega}{Z_{\Lambda_N \setminus A_{k,i},0}^\omega}\\
&\geq \prod_{k=1}^g \prod_{i=1}^{n_k} \frac{c}{\sqrt{\log |t_{k,i}-t_{k,i-1}|}}\,.
\end{align*}
We can now conclude exactly as in Step 2.2. Indeed, the innermost product is of the same type as in \eqref{jensen}, except
that we have an additional log which only helps us.

\end{proof} 

\begin{proof}[Proof of Corollary~\ref{cor_Main}]
The \emph{if} part follows immediately from Theorem~\ref{thm_Main}. To prove the converse, it is enough to bound the indicator function in~\eqref{eq_densityOfVists} by~1.

\end{proof} 

\section*{Acknowledgments} The authors are grateful to M. Lagouge and G. Giacomin for pointing out a mistake in an earlier version of this paper.

\bibliographystyle{plain}
\bibliography{JRV04c}

\end{document}